\documentclass[11pt,a4paper]{amsart}
\usepackage[francais]{babel}
\usepackage[dvips,final]{graphics}
\usepackage{amsmath,amsfonts,amssymb,latexsym}
\usepackage{pstricks}

\usepackage{textcomp}

\newtheorem{thm}{Th\'eoreme}[section]
\newtheorem{cor}[thm]{Corollaire}

\newtheorem{prop}[thm]{Proposition}

\let\ssection=\section
\renewcommand{\section}{\setcounter{equation}{0}\ssection}

\newcommand{\R}{\mathbb{R}}
\newcommand{\Z}{\mathbb{Z}}
\newcommand{\C}{\mathbb{C}}
\newcommand{\N}{\mathbb{N}}

\newcommand{\Q}{\mathbb{Q}}

\newcommand{\B}{\mathcal{B}}
\newcommand{\gn}{\mathfrak{n}}
\newcommand{\gog}{\mathfrak{g}}

\newcommand{\gh}{\mathfrak{h}}

\newcommand{\U}{\mathcal{U}}

\newcommand{\ii}{\textup{\bf{i}}}
\newcommand{\tr}{\textup{Trop }}

\begin{document}

\title[Relevement geometrique de la base canonique...]{Relevement 
geometrique de la base canonique et involution de Sch\"utzenberger}

\bigskip
\bigskip

\author{Sophie Morier-Genoud}
\address{D\'epartement de Math\'ematiques, Universit\'e Claude Bernard 
Lyon I,
69622 Villeurbanne Cedex, France}
\email{morier@igd.univ-lyon1.fr}

\date{}
\maketitle

\begin{abstract}
Soit $G$ un groupe de Lie complexe semisimple simplement connexe, et 
$\B_V$ la base canonique d'un module de Weyl $V$ de $G$. On calcule 
explicitement en terme de para\-m\'etri\-sa\-tion l'action du plus long 
\'el\'ement du groupe de Weyl sur $\B_V$. On utilise pour cela les r\'esultats de 
\cite{berensteinzelevinsky3} sur le rel\`evement g\'eom\'etrique.  
\end{abstract}
\begin{abstract}
Let $G$ be a complex simply connected semisimple Lie group, and let 
$\B_V$ be the canonical base of a Weyl module $V$ of $G$. We calculate 
explicitely the action of the longest element $w_0$ of the Weyl group on 
$\B_V$ in terms of parametrizations. The method is based on results of 
\cite{berensteinzelevinsky3} on the geometric lifting.  
\end{abstract}

\section{Introduction}
Soit $G$ un groupe de Lie complexe semisimple simplement connexe et $B$ 
un sous groupe de Borel. L'alg\`ebre $R$ des fonctions r\'eguli\`eres sur la 
vari\'et\'e de drapeau $G/B$ est engendr\'ee par une base dite base canonique 
duale, \cite{lusztig1}, qui poss\`ede des propri\'et\'es remarquables de 
compatibilit\'es
 avec certaines filtrations. Cette base admet deux syst\`emes de 
para\-m\'etri\-sa\-tions,
 \cite{berensteinzelevinsky3}, la para\-m\'etri\-sa\-tion de Lusztig, et 
la para\-m\'etri\-sa\-tion en cordes.
 On introduit en \ref{defPhi},
 un morphisme $\Phi_\lambda $ qui, \`a un automorphisme de diagramme 
pr\`es, correspond \`a l'action de $w_0$, o\`u $w_0$ est l'\'el\'ement de longueur 
maximale du groupe de Weyl. Ce morphisme $\Phi_\lambda $ pr\'eservant la 
base canonique, on cherche \`a l'exprimer en terme de 
para\-m\'etri\-sa\-tions. Le r\'esultat obtenu g\'en\'eralise \`a tout groupe $G$ et pour tout choix 
d'une d\'ecomposition r\'eduite de $w_0$ le r\'esultat de \cite[\S 
8]{berensteinzelevinsky2}. Un
 r\'esultat remarquable est que pour un choix particulier des 
para\-m\'etri\-sa\-tions, on obtient une expression affine. La combinatoire de la 
base canonique g\'en\'eralise celle des tableaux de Young. Dans cette 
g\'en\'eralisation, le r\'esultat pr\'esent\'e peut-tre vu comme un analogue de 
l'involution de Sch\"utzenberger. La 
m\'ethode utilis\'ee est celle du rel\`evement g\'eom\'etrique,  objet introduit  
et utilis\'e dans une s\'erie d'articles, \cite{berensteinzelevinsky1}, 
\cite{berensteinzelevinsky3}, qui \'etablit un lien entre les combinatoires 
de la base canonique et la g\'eom\'etrie de certaines vari\'et\'es totalement 
positives. A l'aide de ce r\'esultat, on montre dans \cite{CMMG} que l'on 
peut r\'ealiser les c\^ones de Lusztig en termes de param\'etrisations d'une 
partie de la base canonique. Une autre application de ce r\'esultat 
concerne une g\'en\'eralisation de  \cite{SaMajeste} et permet en particulier de 
fournir explicitement des d\'eg\'en\'erescences semi-toriques pour les 
vari\'et\'es de Richardson. Ces r\'esultats sont plus largement d\'evelopp\'es dans  
\cite{moi}. 
 
\section{Notations et pr\'eliminaires}
 \subsection{}
Soit $G$ un groupe de Lie complexe semisimple et simplement connexe. On 
fixe un tore $T$ et  un sous groupe de Borel $B$ 
de $G$. Soit $N$  le radical unipotent de $B$. On note $B^-$ le sous 
groupe de Borel oppos\'e et $N^-$ son radical unipotent. On note 
respectivement $\gog$, $\gh$, $\gn$, $\gn^-$, les $\C$-alg\`ebres de Lie de 
$G$, $T$, $N$, $N^-$.  On a la d\'ecomposition triangulaire 
$\gog=\gn^-\oplus \gh\oplus \gn $. Soit $\lbrace\alpha _i
\rbrace_{1\leq i\leq n}$ une base du syst\`eme de racines correspondant \`a 
cette d\'ecomposition,  $n$ \'etant le rang de 
$\gog$. On note $P$ le r\'eseau des poids engendr\'e par les poids 
fondamentaux $\varpi _i$, $1\leq i\leq n$, et 
$P^+:=\sum_i\N.\varpi _i$ le 
semigroupe des poids entiers dominants. La matrice de Cartan associ\'ee \`a 
$\gog$ est not\'ee $(a_{ij})_{1\leq i,j\leq n}$.
On d\'esigne par $W$ le groupe de Weyl, engendr\'e par les r\'eflexions $s_i$ 
correspondant aux racines simples $\alpha _i$. Un mot r\'eduit pour $w\in 
W$ est une suite d'indices $\ii=(i_1,\cdots, i_l)$ telle
que $w=s_{i_1}\cdots s_{i_l}$, o\`u $\ell (w):=\ell $ d\'esigne la longueur 
de $w$. Soit $w_0$ l'unique \'el\'ement de $W$ de longueur maximale. On 
pose $N:=\ell (w_0)$.

\subsection{}
L'alg\`ebre enveloppante de $\gog$, not\'ee $\U(\gog)$, est engendr\'ee par 
les g\'e\-n\'e\-ra\-teurs $E_i$, $F_i$, $H_i$, $1\leq i\leq n$ et les 
relations de Serre. On note $\B$ la base canonique de $\U(\gn^-)$, et $\tilde 
e_i$, $\tilde f_i$ les op\'erateurs de Kashiwara, \cite{lusztig1}, 
\cite{kashiwara1}. Soit $\lambda $ dans $P^+$, on note $V(\lambda )$ le 
module de Weyl de plus haut poids $\lambda $ et $v_\lambda $ un vecteur de 
plus haut poids de $V(\lambda )$. On note $\B(\lambda ):=V(\lambda )\cap 
\B$ qui est une base, dite canonique, de $V(\lambda)$. On d\'esigne par 
$v^{low}_\lambda $ le vecteur de $V(\lambda )$ de plus bas poids, 
appartenant \`a $\B(\lambda )$. Soit $\ii$  un mot r\'eduit de $w_0$. On adopte 
les notations de \cite[\S 3]{berensteinzelevinsky3}, on note 
$b_\ii:\Z^N_{\geq 0}\rightarrow \B$ la bijection correspondant \`a la 
para\-m\'etri\-sa\-tion de Lusztig de la base canonique et  $c_\ii:\B\rightarrow 
\Z^N_{\geq 0} $ la para\-m\'etri\-sa\-tion en cordes de la base $\B$. On 
 pose $\mathcal{C}_\ii:=c_\ii(\B)$ et $\mathcal{C}_\ii(\lambda 
):=c_\ii(\B(\lambda ))$. Pour deux mots r\'eduits $\ii$ et $\ii'$, on d\'efinit les 
applications de changement de para\-m\'etri\-sa\-tions 
$R_{\ii}^{\ii'}=(b_{\ii'})^{-1}\circ b_\ii: \Z^N_{\geq 0}\rightarrow \Z^N_{\geq 0}$ et
$R_{-\ii}^{-\ii'}=c_{{\ii}'}\circ(c_\ii)^{-1}:{\mathcal C}_{{\bf 
i}}\rightarrow{\mathcal C}_{{\bf i}'}$.

\subsection{} \label{defPhi}
Soit $\omega $ l'automorphisme de $\U(\gog)$ d\'efini sur les 
g\'e\-n\'e\-ra\-teurs par
$\omega (E_i)=F_i$, $\omega (F_i)=E_i$, $\omega (H_i)=-H_i$.
Pour tout $\U(\gog)$-module $V$, soit $V^\omega $ le module tordu par 
l'action de $\omega $, i.e comme espace vectoriel $V\simeq V^\omega $ et
muni de l'action $u\ast v=\omega (u)v$, $u\in \U(\gog)$, $v\in V$. Si 
$V$ est un module simple, alors $V^\omega $ l'est aussi. Ainsi pour tout 
$\lambda \in P^+$, il existe un $\lambda ^\omega \in P^+$ tel que 
$V(\lambda )^\omega $ est isomorphe \`a $V(\lambda^\omega )$. Il existe donc 
un isomorphisme $\Phi _\lambda $ d'espaces vectoriels entre $V(\lambda 
)$ et $V(\lambda^\omega )$ tel que  
$\Phi _\lambda (uv)=\omega(u)\Phi _\lambda(v), u\in \U(\gog), v\in 
V(\lambda )$.
Clairement, $\Phi _\lambda $ envoie un vecteur de plus haut poids sur 
un vecteur de plus bas poids. D'apr\`es le lemme de Schur, ce morphisme 
est unique \`a une constante multiplicative pr\`es. On peut choisir la 
constante pour avoir $\Phi _\lambda (v_\lambda )=v^{low}_{\lambda ^\omega }$.
Par \cite[\S 21]{lusztig1}, on a les propri\'etes suivantes: 

\begin{prop} On a
(i)   $ \lambda ^\omega =-w_0(\lambda )$,   (ii)    $ \Phi _\lambda 
(\B_\lambda )=\B_{\lambda^\omega }$,
(iii)    pour tout $1\leq i\leq n,\Phi _\lambda \tilde {f}_i=\tilde 
{e}_i\Phi _\lambda $
\end{prop}

\section{Rel\`evement g\'eom\'etrique}
Dans cette partie on s'int\'eresse au rel\`evement g\'eom\'etrique de 
l'application $\Phi _\lambda $. En utilisant les r\'esultats de 
\cite{berensteinzelevinsky3}, nous donnons une formule explicite pour l'application 
$b_{\ii}^{-1}\Phi _\lambda c_{\ii}^{-1}$ qui donne la para\-m\'etri\-sa\-tion 
de Lusztig $t'=(t'_1,\cdots, t'_N)$ de l'\'el\'ement $\Phi _\lambda (b)$ en 
fonction de la para\-m\'etri\-sa\-tion en corde $t=(t_1,\cdots ,t_N)$ 
d'un \'el\'ement $b$ de la base canonique. 
\subsection{}
Pour tout $1\leq i\leq n$, on note  $\varphi _i:SL_2\hookrightarrow G$ 
l'injection canonique correspondant \`a la racine simple $\alpha _i$.
Pour $1\leq i\leq n$, on consid\`ere les sous-groupes \`a un param\`etre de 
$G$ d\'efinis par 
$$x_i(t)=\varphi _i\left(
\begin{array}{ccc}
1&t\\
0&1
\end{array}
\right),
\hspace{0.5cm} y_i(t)=\varphi _i\left(
\begin{array}{ccc}
1&0\\
t&1
\end{array}
\right)
,\hspace{0.5cm}t\in \C$$
$$t^{\alpha _i^\vee}=\varphi _i\left(
\begin{array}{ccc}
t&0\\
0&t^{-1}
\end{array}
\right)
,\hspace{0.5cm}t\in \C^*$$Les $x_i(t)$, (resp. $y_i(t)$, $t^{\alpha _i 
^\vee}$) engendrent $N$, (resp.$N^-$, $H$). 
On a les relations de commutation suivantes:
\begin{eqnarray} \label{relcommut}
t^{\alpha _i^\vee}x_j(t')=x_j(t^{a_{ij}}t')t^{\alpha _i^\vee}, 
\hspace{0.3cm}t^{\alpha _i^\vee }y_j(t')=y_j(t^{-a_{ij}}t')t^{\alpha _i^\vee} 
\end{eqnarray}
On d\'efinit deux antiautomorphismes involutifs de $G$, $x\mapsto x^T$, 
appel\'e transposition, et  $x\mapsto x^\iota $, appel\'e inversion, par:
$$
\begin{array}{lll}
x_i(t)^T=y_i(t), & y_i(t)^T=x_i(t),  & (t^{\alpha _i^\vee})^T=t^{\alpha 
_i^\vee}\\[4pt]
x_i(t)^ \iota =x_i(t), & y_i(t)^\iota =y_i(t) ,& (t^{\alpha 
_i^\vee})^\iota =t^{-\alpha _i^\vee}
\end{array}
$$
 
Notons $G_0:=N^-HN$ l'ensemble des \'el\'ements de $G$ qui admettent une 
d\'ecomposition (unique) gaussienne; on \'ecrira $x=[x]_-[x]_0[x]_+$ pour 
$x\in G_0$.

Pour toute suite d'indices  $\ii=(i_1,\cdots, i_m) $ et tout $m$-uplet 
$t=(t_1,\cdots, t_m) $ de $\C^m$, on note:
\vspace{-0.2cm}
$$x_{\ii}(t):=x_{i_1}(t_1)\cdots x_{i_m}(t_m), \hspace{0.3cm} \text{et} 
\hspace{0.3cm} 
x_{-\ii}(t):=y_{i_1}(t_1)t_1^{-\alpha _{i_1} ^{\vee}} \cdots 
y_{i_m}(t_m)t_m^{-\alpha _{i_m} ^\vee} $$
D'apr\`es [2], les $x_\ii$ et $x_{-\ii}$ param\'etrisent des sous vari\'et\'es 
de $G$, et en particulier pour $\ii$ un mot r\'eduit de $w_0$, on a
\begin{thm}
Il existe deux sous vari\'et\'es de $G$, not\'ees $L^{e,w_0}_{>0}$, resp. 
$L^{w_0,e}_{>0}$, telles que pour tout $\ii$ mot r\'eduit de $w_0$, 
l'application $x_{\ii}$, resp. $x_{-\ii}$, r\'ealise une bijection de $\R^N_{>0} 
$ sur  $L^{e,w_0}_{>0}$, resp. $L^{w_0,e}_{>0}$.
\end{thm}

On note $\tilde {R}_\ii^{\ii'}:=x_{\ii'}^{-1}\circ x_\ii$ et $\tilde 
{R}_{-\ii}^{-\ii'}:=x_{-\ii'}^{-1}\circ x_{-\ii}$ les applications de 
changement de para\-m\'etri\-sa\-tion. 
Un r\'esultat important de \cite{berensteinzelevinsky3}  donne que ces 
applications rel\`event g\'eom\'etriquement les applications $R_\ii^{\ii'}$ et 
$R_{-\ii}^{-\ii'}$ de changement de para\-m\'etri\-sa\-tion de la base 
canonique. Plus pr\'ecis\'ement, utilisant les r\'esultats sur les semicorps de 
\cite{berensteinzelevinsky1}, on d\'efinit une transformation de 
"tropicalisation", not\'ee $ [.]_\tr$. Bri\`evement, l'application $ [.]_\tr$ est 
une application entre le semicorps $\Q_{>0}(t_1,\cdots,t_N)$   des 
\textit{expressions rationnelles sans soustraction en }$t_1,\cdots, t_N$ et 
l'ensemble des applications de $\Z^N$ dans $\Z$ et qui consiste \`a 
remplacer la multiplication par l'addition, la division par la soustraction 
et l'addition par l'op\'eration $a\oplus b:=\min (a,b)$. 
On a d'apr\`es \cite{berensteinzelevinsky3},
\begin{thm} \label{ThmTropR}
Les composantes de $(\tilde {R}_\ii^{\ii'})^\vee $ et $(\tilde 
{R}_{-\ii}^{-\ii'})^\vee$ sont des  expressions rationnelles sans soustraction, 
et 
$$ (i)\hspace{0.3cm} [(\tilde {R}_\ii^{\ii'})^\vee  
(t)]_{\tr}=R_\ii^{\ii'}(t) \hspace{0.5cm} (ii) \hspace{0.3cm} [(\tilde 
{R}_{-\ii}^{-\ii'})^\vee  (t)]_{\tr}=R_{-\ii}^{-\ii'}(t) $$

\end{thm}

La notation $(.)^\vee$ signifie que l'on consid\`ere les formules 
analogues dans le dual de Langlands de $G$. On note aussi $[.]_\tr$  
l'application $ [.]_\tr$ appliqu\'ee sur chaque composante.

\subsection{}
Soit $\zeta :L^{w_0,e}_{>0}\rightarrow L^{e,w_0}_{>0} $ d\'efinie par 
$\zeta (x):=[x^{\iota T}]_{+}$.
En utilisant (\ref{relcommut}), on obtient ais\'ement que l'application 
$\zeta $ est bien d\'efinie et plus pr\'ecis\'ement, 
\begin{prop} \label{FormZeta} 
Soit $\ii=(i_1,\cdots, i_N) $ un mot r\'eduit de $w_0$, et $(t_1',\cdots 
,t_N')=(x_\ii^{-1}\circ \zeta  \circ x_{-\ii})(t_1,\cdots ,t_N)$. On a 
alors
$$t'_k=t_k^{-1}\prod _{j>k}t_j^{-a_{i_ji_k}}$$
\end{prop}

On peut maintenant donner la formule du rel\`evement g\'eom\'etrique de $\Phi 
_\lambda $:
\begin{thm}Soit $\ii$ et $\ii'$ deux mots r\'eduits de $w_0$, alors les 
composantes $(x_\ii^{-1}\circ \zeta  \circ x_{-\ii'})^\vee$ sont des 
expressions rationnelles sans soustraction, et
$$ b_{\ii}^{-1}\Phi _\lambda c_{\ii'}^{-1}(t)= [(x_\ii^{-1}\circ \zeta  
\circ x_{-\ii'})^\vee (t)]_\tr+b_{\ii}^{-1}\Phi _\lambda (v_\lambda )$$
\end{thm}

On pose $(l_1,\cdots ,l_N):=b_{\ii}^{-1}\Phi _\lambda (v_\lambda )$. 
Notons que ces constantes peuvent \^etre donn\'ees explicitement, 
\cite{CMMG}. On obtient alors la formule suivante:
\begin{cor}
Pour $(t_1',\cdots ,t_N')=b_{\ii}^{-1}\Phi _\lambda c_{\ii}^{-1}(t_1, 
\cdots,t_N)$,
$$t'_k=l_k-t_k-\sum _{j>k}a_{i_ki_j}t_j$$
\end{cor}

\noindent
\textbf{Preuve:}
On fixe un poids $\lambda $ dans $P^+$. Soit $\Phi 
_{\ii,\ii'}:\mathcal{C} _\ii(\lambda )\rightarrow \Z^N$ une famille d'applications index\'ee 
par deux mots r\'eduits de $w_0$ v\'erifiant les trois conditions 
suivantes:
\begin{enumerate}
\item $\Phi _{\ii,\ii'}(0,\cdots,0)=b_{\ii}^{-1}\Phi _\lambda 
(v_\lambda )$
\item $\Phi _{\ii,\ii'}=R_{\ii''}^\ii\circ \Phi _{\ii'',\ii'}=\Phi 
_{\ii'',\ii'}\circ R_{-\ii'}^{-\ii''}$
\item Pour $\Phi _{\ii,\ii}(t_1,\cdots,t_N)=(t_1',\cdots ,t_N')$, 
$t'_1+t_1$ et  $t'_k$, $k\not=1$,  sont des fonctions de $t_2, \cdots, t_N$.
\end{enumerate}
Le th\'eor\`eme est une cons\'equence de la proposition suivante:

\begin{prop}On a,
\begin{enumerate}
\item[(i)] Si $(\Phi _{\ii,\ii'})$ est une famille v\'erifiant les 
conditions (1), (2), (3), alors
$$ \Phi _{\ii,\ii'}=b_{\ii}^{-1}\Phi _\lambda c_{\ii'}^{-1}$$
\item [(ii)]La famille $(\Phi _{\ii,\ii'})$ d\'efinie par
$$\Phi _{\ii,\ii'}(t)= [(x_\ii^{-1}\circ \zeta  \circ x_{-\ii'})^\vee 
(t)]_\tr+b_{\ii}^{-1}\Phi _\lambda (v_\lambda )$$
v\'erifie les conditions (1), (2), (3).
\end{enumerate}
\end{prop}

\textbf{Preuve:}
Prouvons le point (ii). Remarquant que pour tout $Q$, expression 
rationnelle sans soustraction, $[Q]_\tr(0,...,0)=0$, (1) est clair. Utilisant  
\ref{ThmTropR} et \ref{FormZeta} les points (2) et (3) sont clairs. 
Reste \`a montrer le point (i). Soit $(\Phi _{\ii,\ii'})$  une famille 
v\'erifiant les conditions (1), (2), (3), on peut d\'efinir des applications 
$F_\ii:\B(\lambda )\rightarrow \Z^N_{\geq 0}$ par $F_\ii(b)=\Phi 
_{\ii,\ii'}\circ c_{\ii'}(b)$, $b\in \B(\lambda )$. Ces applications ne 
d\'ependent pas de $\ii'$ d'apr\`es (2). Montrons  par induction sur le poids de 
$b$, que pour tout mot $\ii$, $F_\ii(b)=b_\ii^{-1}(\Phi _\lambda (b))$. 
Si $b=v_\lambda $, c'est clair d'apr\`es (1). Si $b=\tilde {f}_i(b')$, on 
peut choisir un mot r\'eduit $\ii'$ commenant par $i$, 
$\ii'=(i,i'_2,\cdots,i_N')$. On a $F_{\ii'}(b) =\Phi _{\ii',\ii'}\circ c_{\ii'}(\tilde 
{f}_i(b'))=\Phi _{\ii',\ii'}( c_{\ii'}(b')+(1,0,\cdots,0))$. Utilisant 
la condition (3) on a alors $\Phi _{\ii',\ii'}( c_{\ii'}(b')+(1,0,
 \cdots,0)) =\Phi _{\ii',\ii'}\circ c_{\ii'}(b')-(1,0,\cdots,0)$. Ainsi 
$F_{\ii'}(b) = \Phi _{\ii',\ii'}\circ 
c_{\ii'}(b')-(1,0,\cdots,0)=F_{\ii'}(b')-(1,0,\cdots,0)=b_{\ii'}^{-1}(\Phi _\lambda 
(b'))-(1,0,\cdots,0)=b_{\ii'}^{-1}(\tilde {e}_i\Phi _\lambda (b'))=b_{\ii'}^{-1}(\Phi 
_\lambda (\tilde {f}_ib'))=b_{\ii'}^{-1}(\Phi _\lambda (b))$. Puis utilisant 
la condition (2), on a $ F_\ii(b)=b_\ii^{-1}(\Phi _\lambda (b))$, pour 
tout mot r\'eduit $\ii$.

\section{L'involution de Sch\"utzenberger}

\subsection{}

L'involution $i \mapsto i^* $ de l'ensemble $\lbrace 1,\dots,n\rbrace$ 
est d\'efinie par $w_0(\alpha_i)=-\alpha_{i^*}$. Pour une suite d'indice 
$\ii=(i_1,\cdots, i_m) $, on pose $\ii^*:=(i_1^*,\cdots, i_m^*)$. On  
note $\delta$ l'automorphisme d\'efini sur les g\'e\-n\'e\-ra\-teurs par,
$$\delta (E_i)=E_{i^*}, \hspace{0.5cm} \delta (F_i)=F_{i^*}, 
\hspace{0.5cm} \delta (H_i)=H_{i^*}$$
Comme en \ref{defPhi}, pour $\lambda$ dans $P^+$, l'application $\delta 
$ induit un automorphisme de $V(\lambda )$ not\'e $d_\lambda $ qui 
pr\'eserve la base canonique. De plus pour tout \'el\'ement $b_\ii(t)$ de 
$\B(\lambda )$, on a $d_\lambda (b_\ii(t))=b_{\ii^*}(t)$.\\
L'application $\eta _\lambda :=\Phi _\lambda \circ \delta _\lambda $ 
correspond \`a l'involution de Sch\"utzenberger.

\subsection{}
Nous pouvons donner le rel\`evement g\'eom\'etrique de l'application  $\eta 
_\lambda $. Soit $\ii$ un mot r\'eduit de $ w_0$. On pose $(l_1,\cdots 
,l_N):=b_{\ii^*}^{-1}\eta _\lambda (v_\lambda )$.

\begin{cor}
Soit $ (t_1',\cdots ,t_N')=b_{\ii^*}^{-1}\eta _\lambda 
(c_\ii(t_1,\cdots ,t_N))$, alors
$$t'_k=l_k-t_k-\sum _{j>k}a_{i_ki_j}t_j$$
\end{cor}

\end{document}